\documentclass{article}

\usepackage{amsmath}
\usepackage{amssymb}

\newtheorem{thm}{Theorem}
\newtheorem{conj}{Conjecture}
\newtheorem{lem}{Lemma}

\title{Factors of almost squares and lattice points on circles}
\author{Tsz Ho Chan}
\date{}

\begin{document}
\maketitle

\begin{abstract}
In this paper, we consider a conjecture of Erd\H{o}s and Rosenfeld and a conjecture of Ruzsa when the number is an almost square. By the same method, we consider lattice points of a circle close to the $x$-axis with special radii.
\end{abstract}

\section{Introduction and main result}
In [\ref{ER}], Erd\H{o}s and Rosenfeld considered the differences between the divisors of a positive integer $n$. They exhibited infinitely many integers with four ``small" differences and posed the question that any positive integer can have at most a bounded number of ``small" differences. Specifically, they asked
\begin{conj} \label{ques1}
Is there an absolute constant $K$, so that for every $c$, the number of divisors of $n$ between $\sqrt{n} - c \sqrt[4]{n}$ and $\sqrt{n} + c \sqrt[4]{n}$ is at most $K$ for $n > n_0(c)$?
\end{conj}
They also mentioned a conjecture of Ruzsa which is a stronger question:
\begin{conj} \label{ques2}
Given $\epsilon > 0$, is there a constant $K_\epsilon$ such that, for any positive integer $n$, the number of divisors of $n$ between $n^{1/2} - n^{1/2-\epsilon}$ and $n^{1/2} + n^{1/2-\epsilon}$ is at most $K_\epsilon$?
\end{conj}

In [\ref{C}], the author proved the above conjecture of Erd\H{o}s and Rosenfeld [\ref{ER}] for perfect squares and made a tiny progress towards Ruzsa's conjecture for perfect squares:
\begin{thm} \label{thm1}
Any sufficiently large perfect square $n = N^2$ has at most five divisors between $\sqrt{n} - \sqrt[4]{n}(\log n)^{1/7}$ and $\sqrt{n} + \sqrt[4]{n}(\log n)^{1/7}$.
\end{thm}
Here we extend the result slightly to almost squares, namely
\begin{thm} \label{thm2}
Any sufficiently large integer $n$, which can be factored as $(N - a)(N + b)$ for some integers $N$, $a$, $b$ with $0 \le a \le b \le e^{(\log n)^{2/7}}$, has at most eighteen divisors between $\sqrt{n} - \sqrt[4]{n}(\log n)^{1/14}$ and $\sqrt{n} + \sqrt[4]{n}(\log n)^{1/14}$.
\end{thm}
This includes numbers of the form $N^2 - 1$, $N^2 - 4$, $N^2 - N - 6$, ...

\bigskip

Factoring $n = x y$ can also be thought of as finding lattice points on the hyperbola $x y = n$. The above theorems mean that if $x y = n$ has two lattice points within a distance of $(\log n)^{1/7}$ from the point $(\sqrt{n}, \sqrt{n})$, then it can have at most thirty six distinct lattice points within a distance $\sqrt[4]{n}(\log n)^{1/7}$ from the point $(\sqrt{n}, \sqrt{n})$. Similarly, one can consider lattice points on the circle $x^2 + y^2 = n$. It was conjectured that (see [\ref{CG}] for example)
\begin{conj}
For any $\alpha < 1/2$, there exists a constant $C_\alpha$ such that for any $N$ we have
\[
\# \{(a,b): a, b \text{ integers}, a^2 + b^2 = n, N \le |b| < N + n^\alpha \} \le C_\alpha.
\]
\end{conj}
A special case of interest is when $N = 0$:
\begin{equation} \label{circle}
\# \{(a,b): a, b \text{ integers}, a^2 + b^2 = n, |b| < n^\alpha \} \le C_\alpha.
\end{equation}
It is simple to prove (\ref{circle}) for $\alpha \le 1/4$. We extend the range for $|b|$ slightly in special cases of $n$ by showing
\begin{thm} \label{circle1}
For sufficiently large perfect squares $n = N^2$,
\[
\# \{(a,b): a, b \text{ integers}, a^2 + b^2 = n, |b| < n^{1/4} (\log n)^{1/7} \} \le 10.
\]
\end{thm}
\begin{thm} \label{circle2}
For sufficiently large $n$, if $n = a_1^2 + b_1^2$ for some $|b_1| \le e^{(\log n)^{2/7}}$, then
\[
\# \{(a,b): a, b \text{ integers}, a^2 + b^2 = n, |b| < n^{1/4} (\log n)^{1/14} \} \le 36.
\]
\end{thm}

\section{Tools}
The main tool of the proofs is the following result of Turk [\ref{T}] on the size of solutions to simultaneous Pell equations.
\begin{thm} \label{pell}
Let $a$, $b$, $c$, $d$ be squarefree positive integers with $a \neq b$ and $c \neq d$ and let $e$ and $f$ be any integers. If $a f = c e$ then we also assume that $a b c d$ is not a perfect square. Then every positive integer solution of
\[
\left\{ \begin{array}{l}
a x^2 - b y^2 = e \\
c x^2 - d z^2 = f
\end{array} \right.
\]
satisfies
\[
\max(x, y, z) < e^{C \alpha^2 (\log \alpha)^3 \gamma \log \gamma}
\]
where $\alpha = \max(a, b, c, d)$, $\beta = \max(|e|, |f|, 3)$, $\gamma = \max(\alpha \log \alpha, \log \beta)$ and $C$ is a large absolute constant.
\end{thm}
As a consequence, Turk [\ref{T}] proved the following
\begin{thm} \label{thm3}
Suppose $[N, N+K]$, where $K \ge 3$, contains three distinct integers of the form $a_i x_i^2$ with positive integers $a_i$ $x_i$ for $i = 1, 2, 3$. Put $H = \max(a_1, a_2, a_3, 3)$. Then, for some absolute constant $C$,
\[
C H^2 (\log H)^3 (H \log H + \log K) (\log H + \log \log K) > \log N.
\]
\end{thm}
We also need a result on almost squares.
\begin{lem} \label{lem1}
Suppose a positive integer $n$ can be factor as $n = (N - a_1)(N + b_1) = (N - a_2)(N + b_2)$ with $0 \le a_1 < a_2$ and $0 \le b_1 < b_2$. Then $a_2 b_2 \ge N$.
\end{lem}

Proof: From $(N - a_1)(N + b_1) = (N - a_2)(N + b_2)$, we get $(b_1 - a_1) N - a_1 b_1 = (b_2 - a_2) N - a_2 b_2$. Then $[(b_2 - a_2) - (b_1 - a_1)] N = a_2 b_2 - a_1 b_1 > 0$. Hence $N \le [(b_2 - a_2) - (b_1 - a_1)] N \le a_2 b_2$.

\section{Proof of Theorem \ref{thm2}}
We may assume that $n$ is not a perfect square. Suppose $n$ has more than eighteen divisors between $\sqrt{n} - \sqrt[4]{n}(\log n)^{1/7}$ and $\sqrt{n} + \sqrt[4]{n}(\log n)^{1/7}$. Say
\begin{equation} \label{start}
n = (N - a_1)(N + b_1) = (N - a_2)(N + b_2) = ... = (N - a_{10})(N + b_{10})
\end{equation}
with $0 \le a_1 < a_2 < ... < a_{10} \le \sqrt[4]{n}(\log n)^{1/14}$, $0 \le b_1 < b_2 < ... < b_{10} \le \sqrt[4]{n}(\log n)^{1/14}$ and $a_1 \le b_1$ (by taking $N = [\sqrt{n}]$). By Lemma \ref{lem1}, we have $a_i b_i \ge N \gg n^{1/2}$ for $2 \le i \le 10$. Hence
\begin{equation} \label{bound}
a_i, b_i \gg n^{1/4} / (\log n)^{1/14}
\end{equation}
By letting $l_i = b_i - a_i$, we have, from (\ref{start}),
\begin{equation} \label{next}
l_1 N - a_1 (a_1 + l_1) = l_2 N - a_2 (a_2 + l_2) = ... = l_{10} N - a_{10} (a_{10} + l_{10}).
\end{equation}
Observe that $0 \le l_1 < l_2 < ... < l_{10}$ as $a_1 (a_1 + l_1) < a_2 (a_2 + l_2) < ... < a_{10} (a_{10} + l_{10})$.
Note that by the assumption of the theorem, we have $a_1, l_1, a_1 + l_1 \le e^{(\log n)^{2/7}}$. Also $l_i N = a_i b_i \le \sqrt{n} (\log n)^{1/7}$ implies $l_i \le 2 (\log n)^{1/7}$. Subtracting the equations in (\ref{next}), we get
\begin{equation} \label{system1}
\left\{
\begin{array}{l}
(l_2 - l_1) N = a_2 (a_2 + l_2) - a_1 (a_1 + l_1) \\
(l_3 - l_1) N = a_3 (a_3 + l_3) - a_1 (a_1 + l_1) \\
....... \\
(l_{10} - l_1) N = a_{10} (a_{10} + l_{10}) - a_1 (a_1 + l_1).
\end{array}
\right.
\end{equation}
Picking the first three of (\ref{system1}) and eliminating $N$, we have
\[
\left\{
\begin{array}{l}
(l_3 - l_1) a_2 (a_2 + l_2) - (l_3 - l_1) a_1 (a_1 + l_1) = (l_2 - l_1) a_3 (a_3 + l_3) - (l_2 - l_1) a_1 (a_1 + l_1) \\
(l_4 - l_1) a_2 (a_2 + l_2) - (l_4 - l_1) a_1 (a_1 + l_1) = (l_2 - l_1) a_4 (a_4 + l_4) - (l_2 - l_1) a_1 (a_1 + l_1)
\end{array}
\right.
\]
Multiplying everything by $4$, completing the squares and rearranging terms, we have
\begin{equation} \label{system2}
\left\{
\begin{array}{l}
(l_3 - l_1)(2 a_2 + l_2)^2 - (l_2 - l_1)(2 a_3 + l_3)^2 = (l_3 - l_2)(2 a_1 + l_1)^2 + (l_2 - l_1)(l_3 - l_2)(l_1 - l_3) \\
(l_4 - l_1)(2 a_2 + l_2)^2 - (l_2 - l_1)(2 a_4 + l_4)^2 = (l_4 - l_2)(2 a_1 + l_1)^2 + (l_2 - l_1)(l_4 - l_2)(l_1 - l_4)
\end{array}
\right.
\end{equation}
Hence
\begin{equation} \label{system3}
\left\{
\begin{array}{l}
(l_4 - l_1)(l_3 - l_1)(2 a_2 + l_2)^2 - (l_4 - l_1)(l_2 - l_1)(2 a_3 + l_3)^2 = [(l_3 - l_2)(2 a_1 + l_1)^2 + (l_2 - l_1)(l_3 - l_2)(l_1 - l_3)](l_4 - l_1) \\
(l_4 - l_1)(l_3 - l_1)(2 a_2 + l_2)^2 - (l_3 - l_1)(l_2 - l_1)(2 a_4 + l_4)^2 = [(l_4 - l_2)(2 a_1 + l_1)^2 + (l_2 - l_1)(l_4 - l_2)(l_1 - l_4)](l_3 - l_1)
\end{array}
\right.
\end{equation}
Let
\[
(l_4 - l_1)(l_3 - l_1) = s_{34} t_{34}^2, \; \; (l_4 - l_1)(l_2 - l_1) = s_{24} t_{24}^2, \; \; (l_3 - l_1)(l_2 - l_1) = s_{23} t_{23}^2
\]
where $s_{34}$, $s_{24}$, $s_{23}$ are squarefree and less than $4 (\log n)^{2/7}$. Then (\ref{system3}) becomes
\begin{equation} \label{system4}
\left\{
\begin{array}{l}
s_{34} X^2 - s_{24} Y^2 = [(l_3 - l_2)(2 a_1 + l_1)^2 + (l_2 - l_1)(l_3 - l_2)(l_1 - l_3)](l_4 - l_1) \\
s_{34} X^2 - s_{23} Z^2 = [(l_4 - l_2)(2 a_1 + l_1)^2 + (l_2 - l_1)(l_4 - l_2)(l_1 - l_4)](l_3 - l_1)
\end{array}
\right.
\end{equation}
where $X = t_{34} (2 a_2 + l_2)$, $Y = t_{24} (2 a_3 + l_3)$ and $Z = t_{23} (2 a_4 + l_4)$. To apply Theorem \ref{pell}, we need to separate the exceptional cases: 1. $s_{34} = s_{24}$, 2. $s_{34} = s_{23}$, 3. $s_{34} s_{24} s_{34} s_{23}$ is a perfect square.

\bigskip

Case 1: $s_{34} = s_{24}$. Then
\[
s_{34} (X^2 - Y^2) = [(l_3 - l_2)(2 a_1 + l_1)^2 + (l_2 - l_1)(l_3 - l_2)(l_1 - l_3)](l_4 - l_1).
\]
If both sides are not zero, then the absolute value of the left hand side is at least $X + Y \gg n^{1/5}$ by (\ref{bound}) while the right hand side is $O(e^{4 (\log n)^{2/7}})$ which cannot happen if $n$ is sufficiently large. Therefore both sides must be zero. This implies $(2 a_1 + l_1)^2 = (l_3 - l_1)(l_2 - l_1)$. Recall $(l_2 - l_1) N = a_2 (a_2 + l_2) - a_1 (a_1 + l_1)$. Multiplying everything by $4$ and completing the square, we have $(l_2 - l_1) (4 N + l_2 + l_1) = (2 a_2 + l_2)^2 - (2 a_1 + l_1)^2$. Putting in $(2 a_1 + l_1)^2 = (l_3 - l_1)(l_2 - l_1)$ and rearranging terms, we get $(l_2 - l_1) (4 N + l_2 + l_3) = (2 a_2 + l_2)^2$. So in this exceptional case, we have $4 N + l_2 + l_3$ is a positive integer of the form $a x^2$ with $a \le 2 (\log n)^{1/7}$.

\bigskip

Case 2: $s_{34} = s_{23}$. This case is similar to case 1 and we have $4 N + l_2 + l_4$ is of the form $a x^2$ with $a \le 2 (\log n)^{1/7}$.

\bigskip

Case 3: $s_{34} s_{24} s_{34} s_{23}$ is a perfect square. This implies $s_{24} = s_{23}$ since they are squarefree numbers. Subtracting the two equations in (\ref{system4}), we have
\[
s_{23} (Z^2 - Y^2) = (l_4 - l_3)(l_1 - l_2)(2 a_1 + l_1)^2 + (l_2 - l_1)(l_3 - l_1)(l_4 - l_1)(l_4 - l_3).
\]
If both sides are not zero, then the absolute value of the left hand side is at least $Z + Y \gg n^{1/5}$ by (\ref{bound}) while the right hand side is $O(e^{4 (\log n)^{2/7}})$ which cannot happen if $n$ is sufficiently large. Therefore both sides must be zero. This implies $(2 a_1 + l_1)^2 = (l_4 - l_1)(l_3 - l_1)$. Recall $(l_3 - l_1) N = a_3 (a_3 + l_3) - a_1 (a_1 + l_1)$. Multiplying everything by $4$ and completing the square, we have $(l_3 - l_1) (4 N + l_3 + l_1) = (2 a_3 + l_3)^2 - (2 a_1 + l_1)^2$. Putting in $(2 a_1 + l_1)^2 = (l_4 - l_1)(l_3 - l_1)$ and rearranging terms, we get $(l_3 - l_1) (4 N + l_3 + l_4) = (2 a_3 + l_3)^2$. So in this exceptional case, we have $4 N + l_3 + l_4$ is a positive integer of the form $a x^2$ with $a \le 2 (\log n)^{1/7}$.

\bigskip

Therefore, aside from these exceptions, we can apply Theorem \ref{pell} and obtain
\[
n^{1/4} \ll X < e^{8 C (\log n)^{6/7} (\log \log n)^5}
\]
which cannot be true when $n$ is sufficiently large. Hence we have a contradiction or one of the exceptions happens. Similarly we can pick the fourth, fifth and sixth equations in (\ref{system1}), argue in the same way and get a contradiction or $4 N + l_i + l_j$ of the form $a x^2$ with $a \le 2 (\log n)^{1/7}$ for some $5 \le i < j \le 7$. Again we can pick the seventh, eighth and ninth equations in (\ref{system1}) and obtain a contradiction or $4 N + l_p + l_q$ of the form $a x^2$ with $a \le 2 (\log n)^{1/7}$ for some $8 \le i < j \le 10$.

\bigskip

If no contradiction is obtained so far, we have three distinct positive integers $4 N + l_a + l_b$, $4 N + l_i + l_j$, $4 N + l_p + l_q$ of the form $a x^2$ with $a \le 2 (\log n)^{1/7}$ for some $2 \le a < b \le 4$, $5 \le i < j \le 7$ and $8 \le p < q \le 10$. This fits the situation of Theorem \ref{thm3}. Applying it with $H = 2 (\log n)^{1/7}$ and $K = 4 (\log n)^{1/7}$, we get
\[
(\log n)^{2/7} (\log \log n)^3 ((\log n)^{1/7} \log \log n) (\log \log n) \gg \log 4N \gg \log n
\]
which is impossible if $n$ is sufficiently large. With this final contradiction, we prove that $n$ cannot have more than eighteen divisors between $\sqrt{n} - \sqrt[4]{n}(\log n)^{1/7}$ and $\sqrt{n} + \sqrt[4]{n}(\log n)^{1/7}$ and hence Theorem \ref{thm2}.
\section{Proof of Theorem \ref{circle1}}

Suppose
\[
\# \{(a,b): a, b \text{ integers}, a^2 + b^2 = n = N^2, |b| < n^{1/4} (\log n)^{1/7} \} > 10.
\]
Then $N^2 = (N - u_1)^2 + v_1^2 = (N - u_2)^2 + v_2^2 = (N - u_3)^2 + v_3^2$ for some integers $0 < u_1 < u_2 < u_3$ and $0 < v_1 < v_2 < v_3 < n^{1/4} (\log n)^{1/7}$. $N^2 = (N - u_i)^2 + v_i^2$ implies $(N - u_i)^2 > N^2 / 2$ and hence $u_i < N / 2$. It also gives $2 N u_i = u_i^2 + v_i^2$ and $(2 N - u_i) u_i = v_i^2$. These imply $N u_i < n^{1/2} (\log n)^{2/7}$ and $u_i < (\log n)^{2/7}$. Also $v_i^2 > N$ and $v_i > \sqrt{N} = n^{1/4}$. Now write $u_i = s_i t_i^2$ where $s_i$ is squarefree. From above, we have $u_i = s_i t_i^2$ divides $v_i^2$. This implies that $s_i t_i$ divides $v_i$. Write $v_i = s_i t_i w_i$, we have
\[
2 N s_i t_i^2 = (s_i t_i^2)^2 + (s_i t_i w_i)^2 \text{ or } 2 N = s_i t_i^2 + s_i w_i^2 \text{ for } i = 1,2,3.
\]
Note $s_i \le s_i t_i^2 < (\log n)^{2/7}$ and $w_i = v_i / (s_i t_i) > n^{1/4} / (\log n)^{2/7}$. Hence
\begin{equation} \label{circlepell}
\left\{
\begin{array}{l}
s_1 w_1^2 - s_2 w_2^2 = s_2 t_2^2 - s_1 t_1^2 \\
s_1 w_1^2 - s_3 w_3^2 = s_3 t_3^2 - s_1 t_1^2
\end{array}
\right.
\end{equation}
If $s_1 = s_2$, then $w_1^2 - w_2^2 = t_2^2 - t_1^2$. If the left hand side is nonzero, then its absolute value is at least $w_1 + w_2 > n^{1/4} / (\log n)^{2/7}$. However the absolute value of the right hand side is at most $(\log n)^{4/7}$. So we must have $w_1 = w_2$ and $t_1 = t_2$ which forces $u_1 = u_2$ contradicting $u_1 < u_2$. Similarly $s_1 = s_3$ is also impossible. The other exception in applying Theorem \ref{pell} is when $s_1 (s_3 t_3^2 - s_1 t_1^2) = s_1 (s_2 t_2^2 - s_1 t_1^2)$ which implies $u_3 = s_3 t_3^2 = s_2 t_2^2 = u_2$ contradicting $u_2 < u_3$.

\bigskip

Therefore we can apply Theorem \ref{pell} and obtain
\[
\frac{n^{1/4}}{(\log n)^{2/7}} < w_1 < e^{C (\log n)^{4/7} (\log \log n)^3 (\log n)^{2/7} \log \log n}
\]
which is impossible when $n$ is sufficiently large. Consequently we cannot have more than ten pairs of $(a,b)$ and hence Theorem \ref{circle1}.

\section{Proof of Theorem \ref{circle2}}

Suppose
\[
\# \{(a,b): a, b \text{ integers}, a^2 + b^2 = n, |b| < n^{1/4} (\log n)^{1/14} \} > 36.
\]
Then we have $n = a_1^2 + b_1^2 = a_2^2 + b_2^2 = ... = a_{10}^2 + b_{10}^2$ with $a_1 > a_2 > ... > a_{10} > 0$, $0 < b_1 < b_2 < ... < b_{10} < n^{1/4} (\log n)^{1/14}$ and $b_1 < e^{(\log n)^{2/7}}$. Let $a_i = a_1 - l_i$ for $i = 2, 3, ..., 10$. From $a_1^2 + b_1^2 = (a_1 - l_i)^2 + b_i^2$, we have
\begin{equation} \label{circleinitial}
0 < 2 a_1 l_i = b_i^2 - b_1^2 - l_i^2.
\end{equation}
As $a_1 > n^{1/2} / 2$ and $b_i < n^{1/4} (\log n)^{1/14}$, we have $l_i < (\log n)^{1/7}$. By eliminating the terms involving $a_1$ in (\ref{circleinitial}), we have
\[
\left\{
\begin{array}{l}
l_3 b_2^2 - l_2 b_3^2 = (l_3 - l_2) (b_1^2 - l_2 l_3) \\
l_4 b_2^2 - l_2 b_4^2 = (l_4 - l_2) (b_1^2 - l_2 l_4).
\end{array}
\right.
\]
Hence
\[
\left\{
\begin{array}{l}
l_3 l_4 b_2^2 - l_2 l_4 b_3^2 = l_4 (l_3 - l_2) (b_1^2 - l_2 l_3) \\
l_3 l_4 b_2^2 - l_2 l_3 b_4^2 = l_3 (l_4 - l_2) (b_1^2 - l_2 l_4).
\end{array}
\right.
\]
Let $l_3 l_4 = s_{34} t_{34}^2$, $l_2 l_4 = s_{24} t_{24}^2$ and $l_2 l_3 = s_{23} t_{23}^2$ where $s_{34}$, $s_{24}$ and $s_{23}$ are squarefree and less than $(\log n)^{2/7}$. Therefore
\begin{equation} \label{circlesystem}
\left\{
\begin{array}{l}
s_{34} X^2 - s_{24} Y^2 = l_4 (l_3 - l_2) (b_1^2 - l_2 l_3) \\
s_{34} X^2 - s_{23} Z^2 = l_3 (l_4 - l_2) (b_1^2 - l_2 l_4)
\end{array}
\right.
\end{equation}
where $X = t_{34} b_2$, $Y = t_{24} b_3$ and $Z = t_{23} b_4$. From (\ref{circleinitial}), we have $b_i^2 > 2 a_1 > n^{1/2}$. Hence $X, Y, Z > n^{1/4}$. Similar to the proof of Theorem \ref{thm2}, we look at the exceptions in applying Theorem \ref{pell}. If $s_{34} = s_{24}$, then
\[
s_{34}(X^2 - Y^2) = l_4 (l_3 - l_2) (b_1^2 - l_2 l_3).
\]
If the left hand side is nonzero, then its absolute value is at least $X + Y > n^{1/4}$. However the right hand side is at most $e^{4 (\log n)^{2/7}}$ which is a contradiction. So we must have the two sides are zero and $b_1^2 = l_2 l_3$. Putting this into (\ref{circleinitial}), we have
\[
2 a_1 l_2 = b_2^2 - l_2 l_3 - l_2^2 \text{ or } (2 a_1 + l_2 + l_3) l_2 = b_2^2.
\]
So the number $2 a_1 + l_2 + l_3$ is of the form $a x^2$ with $a < (\log n)^{1/7}$.

\bigskip

Similarly if $s_{34} = s_{23}$, then we have $2 a_1 + l_2 + l_4$ is of the form $a x^2$ with $a < (\log n)^{1/7}$.

\bigskip

Finally if $s_{34} s_{24} s_{34} s_{23}$ is a perfect square, we must have $s_{24} = s_{23}$ as they are squarefree. Subtracting the two equations in (\ref{circlesystem}), we have
\[
s_{23} Z^2 - s_{24} Y^2 = l_2 (l_3 - l_4) (b_1^2 - l_3 l_4).
\]
Using the same argument as above, we produce $2 a_1 + l_3 + l_4$ is of the form $a x^2$ with $a < (\log n)^{1/7}$.

\bigskip

Aside from these exceptions, we can apply Theorem \ref{pell} to (\ref{circlesystem}) and get
\[
n^{1/4} < X < e^{C (\log n)^{4/7} (\log \log n)^3 (\log n)^{2/7} \log \log n}
\]
which is a contradiction. Therefore we have a contradiction or one of these exceptions happen. Now we repeat the same argument for $a_1^2 + b_1^2 = a_5^2 + b_5^2 = a_6^2 + b_6^2 = a_7^2 + b_7^2$ and $a_1^2 + b_1^2 = a_8^2 + b_8^2 = a_9^2 + b_9^2 = a_{10}^2 + b_{10}^2$ and get $2 a_1 + l_a + l_b$, $2 a_1 + l_i + l_j$ and $2 a_1 + l_p + l_q$ all of the form $a x^2$ with $a < (\log n)^{1/7}$ for some $2 \le a < b \le 4$, $5 \le i < j \le 7$ and $8 \le p < q \le 10$ if we have not already got a contradiction. Then we can apply Theorem \ref{thm3} with $H < (\log n)^{1/7}$ and $K < 2 (\log n)^{2/7}$, and get
\[
C (\log n)^{2/7} (\log \log n)^3 ((\log n)^{1/7} \log \log n) (\log \log n) > \log a_1 \gg \log n.
\]
This is impossible if $n$ is sufficiently large. With this final contradiction, we prove Theorem \ref{circle2}.

Department of Arts and Sciences \\
Victory University \\
255 N. Highland St., \\
Memphis, TN 38111 \\
U.S.A. \\
thchan@victory.edu

\end{document}